\documentstyle[12pt]{article}
\newcommand{\sun}{ \sum_{i=1}^{n}}
\newcommand{\h}{\cal H}
\newcommand{\hn}{{\h}_n}
\newcommand{\ltn}{l^2( {\bf N})}
\newcommand{\fti}{ \{f_i \}_{i \in I}}
\newcommand{\gti}{ \{g_i \}_{i \in I}}
\newcommand{\ftu}{ \{f_i \}_{i=1}^{\infty}}
\newcommand{\fn}{ \{f_i \}_{i=1}^n }
\newcommand{\fo}{ \{f_i \}_{i=n+1}^{\infty}}
\newcommand{\fno}{ \{f_i \}_{i=m+1}^{\infty}}
\newcommand{\gno}{ \{g_i \}_{i=m+1}^{\infty}}
\newcommand{\siu}{ \sum_{i=1}^{\infty}}
\newcommand{\gik}{ \{ g_{i_k} \}_{k=1}^{m-n}}
\newcommand{\gtu}{ \{g_i \}_{i=1}^{\infty}}
\newcommand{\etu}{ \{ e_i \}_{i=1}^{\infty}}
\newcommand{\snu}{ \sum_{n=1}^{\infty} }
\newcommand{\smo}{ \sum_{i=m+1}^{\infty}}
\newcommand{\sui}{ \sum_{i \in I}}
\newcommand{\sjn}{\sum_{j=1}^{n}}
\newcommand{\suin}{\sum_{i=1}^{n}}
\newcommand{\sne}{\sum_{i=n+1}^{\infty}}
\newcommand{\sjne}{ \sum_{j=1}^{n+1}}
\newcommand{\sju}{ \sum_{j=1}^{\infty}}
\newcommand{\nf}{||f||}
\newcommand{\nft}{ ||f||^2}
\newcommand{\ffj}{ <f, f_j> }
\newcommand{\ffi}{ <f, f_i> }
\newcommand{\frn}{ \sjne | \ffj |^2 }
\newcommand{\fni}{ f_i^n}
\newcommand{\ftni}{ \{ \fni \}_{i=1, n=1}^{n+1, \infty} }

\newcommand{\ftns}{ \{f_i \}_{i \in N - \sigma}}
\newcommand{\ts}{T^{*}}
\newcommand{\us}{U^{*}}
\newcommand{\cin}{ c_i^n}
\newcommand{\ctu}{ \{c_i\}_{i=1}^{\infty}}
\newcommand{\xtu}{ \{x_i \}}
\newcommand{\ytu}{ \{y_i \}}
\newcommand{\ztu}{ \{z_j \}}
\newcommand{\si}{S^{-1}}
 
\begin{document}
\title{ Frames containing a Riesz basis and preservation of this property
under perturbations.} 
\author{ Peter G. Casazza and Ole Christensen \thanks{ The first author
acknowledges support by NSF Grant DMS-9201357 and a grant from the Danish 
Research Foundation. 1993 Mathematics Subject Classification:
Primary 42C99, 46C99}}
\maketitle
\begin{abstract} Aldroubi has shown how one can construct any frame $\gtu$ 
starting with one frame $\ftu $,using a bounded operator $U$ on $l^2(N)$. 
We study the overcompleteness of the frames in terms of properties of $U$.
We also discuss perturbation of frames in the sense that two frames are
``close'' if a certain operator is compact. In this way we obtain an 
equivalence relation with the property that members of the same 
equivalence class have the same overcompleteness. On the other hand we show
that perturbation in the Paley-Wiener sense does not have this property. \\
Finally we construct a frame which is norm-bounded below, but which does not
contain a Riesz basis.The construction is based on the delicate difference
between the unconditional convergence of the frame representation, and the
fact that
a convergent series in the frame elements need not  converge unconditionally. 
\end{abstract} 
\section{Introduction.} The introduction of frames for a Hilbert space $\h$
goes back to the paper \lbrack DS\rbrack \ from 1952, where they are used
in nonharmonic Fourier analysis. A frame is a family of elements in $\h$
which can be considered as an ``overcomplete basis'': every element in $\h$
can be written as a linear combination of the frame elements, with square
integrable coefficients, which do not need to be unique. A natural 
theoretical question (which is also
important for applications, e.g., representation of an operator using a basis)
is how far frames are away from bases, i.e., one may ask questions like
\\ \\ 1) does a frame contain a Riesz basis?  \\ 
2) which conditions imply that a frame just consists of a Riesz basis
plus finitely many elements ? \\ 
3) what happens with the overcompleteness if the frame elements are
perturbed? \\ \\
The reason for the interest in Riesz bases and not just bases is that
frames and Riesz bases are closely related: a 
Riesz bases is just a frame, where the elements are $\omega$-independent.
\\ \\ 
Some answers has been found by Holub \lbrack H\rbrack \ , who concentrates
on the second question. Here we go one step further, in that we are mainly
interested in frames which just contain a Riesz basis. For such frames
one defines the excess as the number of elements one should take away to
obtain a Riesz basis. \\
In the first part of the paper we apply a result of Albroubi \lbrack A\rbrack
\ , explaining how one can map a frame onto  another 
using a bounded operator $U$ on $l^2$. Our results concern the relation
between the frames involved and properties of $U$. Independent of that we
construct a norm-bounded frame not containing a Riesz basis. \\ \\
In section 3 we concentrate on the third question. We introduce the concept
``compact perturbation''. This leads to an equivalence relation on the set
of frames, with the property that frames in the same equivalence class have 
the same overcompleteness properties; this means, that if a frame contains
a Riesz basis then all members in the class contain a Riesz basis, and all
those frames have the same excess. \\
Finally we show that perturbation in the Paley-Wiener sense \lbrack C3\rbrack
\  not has this plesant property. 
\section{Frames containing a Riesz basis.}
Let $\h$ be a separable Hilbert space. A family $\fti$
is called a {\it frame} for $\h $ if
$$ \exists A,B >0: \ A \nft \leq \sui | <f,f_i>|^2 \leq B \nft , \forall 
f \in \h .$$ A and B are called {\it frame bounds}. \\ \\
A {\it Riesz basis} is a family of elements which is the image of an 
orthonormal
basis by a bounded invertible operator. For families with the natural
numbers as index set there is an equivalent characterization \lbrack Y\rbrack:
$\ftu$ is a Riesz basis if there exist numbers $A,B >0$ such that
$$ (1) \ \ \ \ A \sum_{i=1}^n |c_i |^2 \leq || \sum_{i=1}^n c_i f_i ||^2 \leq
B \sum_{i=1}^n |c_i |^2 ,$$
for all finite sequences $c_1, ...c_n .$ \\ Also, a basis $\ftu$ is a Riesz
basis if and only if it is unconditional (meaning that if  
$\siu c_i f_i $ converges for some coefficients $\{c_i \}$, 
then it actually converges unconditionally) and 
$ 0 < inf_i ||f_i || \leq sup_i ||f_i || < \infty .$ \\ \\
There is a close connection between frames and Riesz bases:
$$ \ftu \ \mbox{is a Riesz basis} \ \Leftrightarrow [ \ftu \ \mbox{ is a frame
and } \ \sum_{i=1}^{\infty} c_i f_i =0 \Rightarrow c_i =0, \forall i .]$$
If $\ftu$ is a Riesz basis, then the numbers $A,B$ appearing in (1) are 
actually frame bounds. \\ \\   
If $\fti$ is a frame (or if only the upper frame condition is satisfied)
then we define the {\it pre-frame operator} by
$$T : l^2(I) \rightarrow \h , \ \ T \{c_i \}:= \sui c_i f_i .$$
The operator $T$ is bounded. Composing $T$ with its adjoint 
$$T^* : \h \rightarrow l^2 (I), \ \  T^*f= \{<f, f_i > \}_{i \in I} $$
we get the {\it frame operator} 
 $$S=TT^* : \h \rightarrow \h , \ \  Sf := \sui <f,f_i> f_i ,$$
which is a bounded and invertible operator. 
 This immediately leads to the {\it frame decomposition} ; every $f \in \h$
can be written as 
$$f= \sui <f, S^{-1}f_i >f_i ,$$ where the series converges unconditionally.
So a frame has a property similar to a basis: every element in $\h$ can be
written as a linear combination of the frame elements. For more 
information about basic properties of frames we refer to the original paper 
\lbrack DS\rbrack \ and the research tutorial \lbrack HW\rbrack. \\ \\ 
The main difference between a frame $\ftu$ and a basis 
is that a frame can be overcomplete, so it might happen that
$f \in \h$ has a representation $f = \siu c_i f_i $ for some coefficients $c_i$
which are different from the {\it frame coefficients} $ <f, S^{-1}f_i> $.
In applications one might wish not to have ``too much redundancy''. In that
spirit Holub \lbrack H\rbrack \  discusses {\it near-Riesz bases}, i.e. frames
$\ftu$ consisting of a Riesz basis $\{f_i \}_{i \in N- \sigma}$
plus finitely many elements $\{f_i \}_{i \in \sigma}$. The number of
elements in $\sigma$ is called the {\it excess}. 
Let us denote the kernel of the operator $T$ by $N_T$.  
If $\ftu$ is a frame, then
$$ \ftu \ \mbox{is a near-Riesz basis} \Leftrightarrow N_T \ \mbox{has finite
dimension}$$\ $$ \Leftrightarrow \ftu  \mbox{ is unconditional} .$$
The first of the above biimplications is due to Holub \lbrack H\rbrack,
who also proves the second under the assumption that the frame is 
norm-bounded below. The generalization above is proved by the authors in
\lbrack CC\rbrack. \\ \\ If the conditions above are satisfied, 
then the excess is equal to dim$(N_T)$. \\ \\ If $\mbox{dim}(N_T) = \infty$, 
two things can happen: $\ftu$ consists of a Riesz
basis plus infinitely many elements (in which case we will say that $\ftu$
has infinite excess) or $\ftu$ does not contains a Riesz basis at all.
In the present paper we concentrate on frames which contain a Riesz basis.
Every frame can be mapped onto such a frame (in  fact, onto an arbitrary
frame) using a construction of Aldroubi \lbrack A\rbrack \ , which we 
now shortly describe. \\ \\
Let $\ftu$ be a frame and $U: \ltn \rightarrow \ltn$ a bounded operator.
Let $\{u_{i,j}\}_{i,j \in N}$ be the matrix for $U$ with respect to some 
basis. Define the family $\gtu \in \h$ by  
$$ g_i = \sju u_{i,j} f_j .$$
By an abuse of notation we will sometimes write $\gtu = U\ftu$.
A result of Aldroubi (differently formulated) states that
$$\gtu \ \mbox{is a frame} \Leftrightarrow \exists \gamma >0 :
 || U \ts f || \geq \gamma \cdot || \ts f || , \ \ \forall f \in \h . $$ 
This is not too complicated: the boundedness of $U$ implies that 
$\gtu$ satisfies the upper frame condition, and the condition above is
just a different expression for the lower condition. But it is important
that {\it every} frame $\gtu$ can be generated in this way, i.e., given the
frame $\gtu$ we just have to find the operator $U$ mapping $\ftu$ to $\gtu.$
\\ \\ In connection with Aldroubis construction there are (at least) two
natural questions related to Holubs work: how is the excess of $\gtu$
related to that of $\ftu$, and which conditions imply that $\gtu$ actually
is a Riesz basis? We shall give answers to both questions in 
this section. \\ \ 
The definition of $\gtu$ immediately shows that
$$ \{ <g_i, f >\} = U \{ <f_i,f> \} , \ \ \forall f \in \h ;$$
this is true whether or not $\gtu$ builds a frame.
The formula leads  to an expression for the pre-frame operator associated 
with $\gtu$. We let $U^{T}$ denote the transpose of $U$ and $\overline{U}$
be the operator corresponding to the matrix where all entries in the matrix 
of $U$ are complex
conjugated. Then, given $f \in \h , \ctu \in \ltn $, we have
$$ < \siu c_i g_i ,f> = \siu c_i <g_i,f> = \siu c_i \overline{<f, g_i>}$$
\ $$=  \siu < \ctu ,\overline{U} \{<f,f_i> \} > 
= <\ctu ,\overline{U} \ts f> = < T U^{T} \ctu, f > .$$  
It follows that
$$ \siu c_i g_i = TU^{T} \ctu , \ \ \forall \ctu \in \ltn .$$
So if $\gtu$ contains a Riesz basis, then its excess is equal to 
$dim( N_{TU^T})$. For the calculation of this number we need a lemma.
Corresponding to an operator $V$ we denote its range by $R_V$. \\ \\
{\bf Lemma 1:} {\it Let $X,Y$ be vector spaces and $V : X \rightarrow Y$ 
a linear mapping. Given a subspace  $Z \subseteq Y $, define
$$ V^{-1}(Z) := \{x \in X \ | \ Vx \in Z \}.$$  Then
$$ \mbox{dim}(V^{-1}(Z))= \mbox{dim}(Z \cap R_V ) + \mbox{dim}(N_V) .$$}
{\bf Proof:} Let $ \ytu \subseteq Y $ be a basis for $Z \cap R_V $ and
take $\xtu \subseteq X$ such that $Vx_i =y_i .$ Now, if $x \in X $ and 
$Vx \in Z$ then we can find coefficients $\{c_i \}$ such that
$Vx = \sum c_i y_i = V \sum c_i x_i $, i.e., 
$x \in \overline{\mbox{span}} \xtu +
N_V .$ Let now $\ztu$ be a basis for $N_V $. Corresponding to our 
element $x \in X $ with $Vx \in Z $ we can now also find coefficients
$\{d_j \} $ such that $x = \sum c_i x_i + \sum d_j  z_j .$ So the 
independent set $ \xtu \cup \ztu $ spans $ V^{-1}(Z) .$ \\ \\
{\bf Theorem 2:} {\it $\mbox{dim} (N_{TU^{T}})= \mbox{dim} (R_{U^T} \cap N_T) +
\mbox{dim}({R}^{\bot}_U) .$} \\ \\
{\bf Proof:}  
$$  \{ \ctu \ | \ TU^T \ctu =0 \} 
= \{ \ctu \ | \ U^T \ctu \in N_T \} = (U^T)^{-1} (N_T) .$$ 
Now the result follows from Lemma 1 and the observation
$$ \mbox{dim}(N_{U^T})= \mbox{dim}(N_{\us})= \mbox{dim}({R_U}^{\bot}) .$$  \\
So if $ \gtu$ actually is a frame containing a Riesz basis, then Theorem 2
gives a recipe for calculation of the excess. In particular, if $\ftu$ is a 
near-Riesz basis and $R_U$ has finite codimension, then also $\gtu$ is a
near-Riesz basis. \\ \\
{\bf Proposition 3:}
{\it $$ \gtu \ \mbox{ is a Riesz basis} \Leftrightarrow 
\overline{U}: R_{\ts} \rightarrow
\ltn \ \mbox{is surjective} .$$}
{\bf Proof:}
The standing assumption ``$U$ bounded'' implies that $\gtu$ satisfies the
upper condition, so $ \gtu$ is a Riesz basis if and only if 
$$ \exists m >0 : \ \ \sun |c_i|^2 \leq m \cdot  || \sun c_i g_i ||^2 $$
for all finite sequences $c_1,...c_n .$ By \lbrack Y, p.155\rbrack \ this
condition is satisfied if and only if the moment problem
$$ < f, g_i> = c_i , \ \ i \in N$$ has at least one solution $f$ whenever 
$\ctu \in \ltn $ ,i.e., if and only if 
$\overline{U} \{ \ffi \} = \ctu $ has at least one
solution whenever $\ctu \in \ltn .$ The last condition means exactly 
that $\overline{U}$ is surjective
considered as a mapping from the subspace $R_{\ts}$ of $\ltn$ onto $\ltn $.  
\\ \\
More generally one may wish that the frame at least contain a Riesz basis.
As shown in \lbrack C2\rbrack \ it is the case for a {\it Riesz frame}, which
is a frame with the property that every subfamily is a frame for its closed 
linear span, with a common lower bound. \\ 
It is easy to construct a frame which does not contains a Riesz basis if
one allows a subsequence of the frame elements to converge against 0 in norm.
We now present an example showing that the same can be the case for a frame
which is norm-bounded below. Our approach is complementary to recent work
of Seip \lbrack Se\rbrack, who proves that there exist frames of complex 
exponentials for $L^2( - \pi , \pi )$ which do not contain a Riesz basis.
While Seip relies on the theory for sampling and interpolation our approach 
is more elementary, just using functional analysis. Furthermore our
construction puts focus on a different point, namely the difference between
convergence and unconditional convergence of an expansion in the 
frame elements. \\ 
\\ {\bf Proposition 4:} {\it There exists a frame for $\h$ made up of norm one 
vectors , which has no subsequence which is a Riesz basis.} \\ \\
The proof needs several lemmas, so let us shortly explain the basic idea. As 
we have seen, $\sui c_i f_i $ converges unconditionally for every set of
frame coefficients $\{c_i \} $. But nothing guarantes that convergence of
$ \sui c_i f_i $ implies unconditional convergence for general coefficients
$\{c_i\} $. Our proof consists in a construction of a frame where no total
subset is unconditional, and hence not a Riesz basis. Technically the first 
step is  to decompose $\h$ 
into a direct sum of finite dimensional subspaces of increasing dimension.
The idea behind the proof might be useful in other situations as well.
\\ \\ {\bf Lemma 5:} {\it Let $\{e_i \}_{i=1}^{n}$ be an orthonormal 
basis for a finite dimensional space $\hn$. Define
$$f_j = e_j - \frac{1}{n} \sun e_i \ \ \mbox{for} \ j=1,..n $$   \ 
$$f_{n+1}= \frac{1}{\sqrt{n}} \sun e_i .$$
Then
 $$  \frn = \nft , \forall f \in 
\hn .$$} \\ \\ {\bf Proof:}  
Given $f \in \hn $, write $f= \suin a_i e_i , \ a_i = <f, e_i > .$ If we let
$P$ denote the orthogonal projection onto the unit vector
$\frac{1}{\sqrt{n}} \suin e_i $, then
$$ Pf = \frac{1}{n} <f, \suin e_i > \suin e_i = \frac{\suin a_i}{\sqrt{n}}
\frac{1}{\sqrt{n}} \suin e_i .$$ Therefore
$$ || Pf||  = \frac {| \suin a_i |^2}{n} = | <f, f_{n+1} >|^2 .$$
Also 
$$ || (I-P)f||^2 = ||f-Pf||^2 = || \suin a_i e_i - \frac{\sjn a_j}{n}
\suin e_i ||^2 $$ \ $$ = || \suin (a_i - \frac{\sjn a_j}{n}) e_i ||^2 =
\suin | a_i - \frac{\sjn a_j}{n} |^2 = \suin |<f, f_i >|^2 .$$
Putting the two results together we obtain
$$ ||f||^2 = ||Pf||^2 + || (I-P)f||^2 = \sum_{i=1}^{n+1} |<f, f_i >|^2 $$
and the proof is complete.  
 \\ \\  
Given a sequence $\{ g_i \}_{i \in I} \subseteq \h $ its {\it unconditional
basis constant} is defined as the number
$$ \mbox{sup} \{ || \sui {\sigma}_i c_i g_i || \ | 
\ || \sui c_i g_i || =1  \ \mbox{and} \ {\sigma}_i = \pm 1, \forall i \} .$$
As shown in \lbrack Si\rbrack , a total family $\{g_i\}_{i \in I}$ consisting
of non-zero elements is an unconditional basis for $\h$ if and only if it has 
finite unconditional basis constant. \\ \\     
{\bf Lemma 6:} {\it Define $\{f_1, ....f_{n+1} \}$ as in Lemma 5. 
Any subset of $\{f_1, f_2,.. f_{n+1} \} $ which spans $\hn$
has unconditional basis constant greater than or equal to $\sqrt{n-1} -1 $.}
\\ \\ {\bf Proof:} Since $ \sun f_i =0$, any subset of $\{f_1,..f_{n+1} \}$
which spans $\hn$ must contain $n-1$ elements from $\{f_1,...f_n \}$ plus 
$f_{n+1}.$ By the symmetric construction it is enough to consider the family
$\{f_1,..f_{n-1},f_{n+1} \} .$ We have \\ \\
$ || \sum_{i=1}^{n-1} f_i || = 
||\sum_{i=1}^{n-1} e_i - \frac{n-1}{n} \sun e_i || =
|| (1- \frac{n-1}{n} \sum_{i=1}^{n-1} e_i - \frac{n-1}{n} e_n ||  \\ \\
= || \frac{1}{n} \sum_{i=1}^{n-1} e_i - \frac{n-1}{n}e_n || =
\sqrt{ \frac{n-1}{n^2} + \frac{(n-1)^2}{n^2}} = \frac{1}{n} \sqrt{n(n-1)}
\leq 1 . $ \\ \\
Now consider $|| \sum_{i=1}^{n-1} (-1)^n f_i || $; if $n$ is odd this number
is equal to $ || \sum_{i=1}^{n-1} (-1)^n e_i || = \sqrt{n-1}$, and if $n$ is
even it is equal to 
$$ || \sum_{i=1}^{n-1} (-1)^i e_i - \frac{1}{n} \sun e_i || \geq
|| \sum_{i=1}^{n-1} (-1)^i e_i || - || \frac{1}{n} \sun e_i || \geq 
\sqrt{n-1} - \frac{\sqrt{n}}{n} \geq \sqrt{n-1}- 1.$$ That is, in all cases,
$$ || \sum_{i=1}^{n-1} (-1)^n f_i || \geq \sqrt{n-1}-1 .$$
Combining this with the norm estimate $|| \sum_{i=1}^{n-1} f_i || \leq
1 $ \ it follows that the unconditional basis constant of 
$\{f_1,...f_{n-1} \}$ is greater than or equal to $\sqrt{n-1}-1$, 
so clearly the same is true for $\{f_1,...f_{n-1},f_{n+1} \} $. \\ \\
Now we are ready to do the construction for Proposition 4. Let $\etu$ be an
orthonormal basis for $\h$ and define
 $$ \hn := \mbox{span} \{e_{\frac{(n-1)n}{2}+1}, e_{\frac{(n-1)n}{2}+2},...
e_{\frac{(n-1)n}{2}+ n} \} .$$
So ${\h}_1 = \mbox{span} \{e_1 \}, \ \ 
{\h}_2 = \mbox{span} \{e_2, e_3 \} , {\h}_3 = \mbox{span} \{e_4, e_5, e_6 \},..
... $. \\ By construction, $$ \h =  (\snu \bigoplus \hn )_{\h} .$$
That is, $g \in \h \Leftrightarrow g = \snu g_n , g_n \in \hn $, \ and 
$ ||g||^2 = \snu ||g_n ||^2$. We refer to \lbrack LT\rbrack \ for details
about such decompositions.  \\ \\
For each space $\hn$ we construct the sequence $\{f_i^n\}_{i=1}^{n+1} $ as in
Lemma 5, starting with the orthonormal basis $ \{e_{\frac{(n-1)n}{2}+1},...
e_{\frac{(n-1)n}{2}+n} \}.$ Specifically, given $n \in N $,
$$ f_i^n = e_{\frac{(n-1)n}{2}+i} - \frac{1}{n} \sjn e_{\frac{(n-1)n}{2}+j} 
\ , \ \ 1\leq i \leq n $$ \
$$ f_{n+1}^n = \frac{1}{\sqrt{n}} \sjn e_{\frac{(n-1)n}{2} +j} .$$
{\bf Lemma 7:} {\it $\ftni$ is a frame for $\h$, with bounds
$ A=B=1$.} \\ \\
{\bf Proof:} Write $g \in \h $ as $ g = \snu g_n, \ g_n \in \hn .$ Given 
$n \in N$ it is clear that 
$$ <g, \fni > = <g_n, \fni > \ \mbox{for} \ i=1,...n+1 .  $$
>From this calculation it follows that 
$$ \snu \sum_{i=1}^{n+1} |<g, \fni >|^2 = 
\snu \sum_{i=1}^{n+1} | <g_n, \fni >|^2 = \snu  ||g_n||^2 =  ||g||^2, $$
where we have used Lemma 5.  \\ \\
{\bf Lemma 8:} {\it No subsequence of $\ftni$ is a Riesz basis for $\h$.} \\ \\
{\bf Proof:} Any subsequence of $\ftni$ which spans $\h$ must contain 
$n$ elements from
$ \{\fni \}_{i=1}^{n+1} $ and so by Lemma 6, its unconditional basis constant
is greater than or equal to $\sqrt{n-1} -1 $ for every $n$. That is, the 
unconditional basis constant is infinite,  hence the subsequence can not  be 
an unconditional basis for $\h$. \\ \\
Lemma 7 and Lemma 8 proves Proposition 4.
It would be interesting to determine whether Proposition 4 still holds 
if one only considers classes of frames with a special structure, for
example Weyl-Heisenberg frames, wavelet frames, or frames consisting of
translates of a single function. \\ \\
{\bf Remark:} The {projection method} developed in \lbrack C1, C2\rbrack \
can be used to calculate the frame coefficients if a certain technical
condition is satisfied. The ``block structure'' of the frame $ \ftni$
constructed here shows that the projection method can be used. As shown
in \lbrack C2\rbrack \ the method can also be used for every Riesz frame.
The two questions, i.e. the question about containment of a Riesz basis and
the question whether the projection method works, don't seem to be strongly
related. 
\section{Excess preserving perturbation.}
At several places in the following we need results for perturbation of
frames and Riesz bases. We denote the frames by $\ftu , \gtu $, usually
with the convention that $\ftu$ is the frame we begin with, and $\gtu$ is
the perturbed family.  
Common for all the result is that they can be formulated using
the {\it perturbation operator} $K$ mapping a sequence $\{ c_i \}$ of
numbers to $\sum c_i (f_i - g_i) $. \\ \\
{\bf Theorem 9:} {\it Let $\ftu , \gtu \subseteq \h .$ \\ \\
a) If $\ftu $ is a frame for $\h$ and $K$ is compact as an operator from 
$l^2(N)$ into $\h$, then $\gtu$ is a frame for its closed linear span. \\ \\
b) Suppose $\ftu$ is a frame for $\h$ with bounds $A,B$. If there exist
numbers $\lambda , \mu \geq 0 $ such that $\lambda + \frac{\mu}{\sqrt{A}} 
< 1$ and $$|| \sum c_i (f_i - g_i ) || \leq \lambda \cdot || \sum c_i f_i ||
+ \mu \sqrt{\sum |c_i |^2 }$$ for all finite sequences $\{c_i \}$, then $\gtu$
is a frame for $\h$ with bounds \\ $A(1-(\lambda + \frac{\mu}{\sqrt{A}}))^2,
B(1+ \lambda + \frac{\mu}{\sqrt{B}})^2.$} \\ \\ For the proofs we refer 
to \lbrack C3, CH\rbrack. \
Several variations are possible. If in b) we just assume that $\ftu$
is a Riesz basis for $\overline{span} \ftu$, then $\gtu$ is also a 
Riesz basis for its closed linear span. Also observe, that if $\ftu$ 
is a frame and $\sigma \subseteq N $ is finite, then $\ftns$ is a frame for 
$\overline{span} \ftns$; this is a consequence of a). The next result 
connects Theorem 9 with the question about overcompleteness. \\ \\
{\bf Theorem 10:} {\it Suppose that $\ftu$ is a frame containing a Riesz
basis, that $\gtu$ is total, and that $K$ is compact as a mapping from
$l^2(N)$ into $\h$. Then $\gtu$ is a frame for $\h$ containing a Riesz basis,
and the frames $\ftu$ and $ \gtu$ have the same excess.} \\ \\
{\bf Proof:} First assume that $\ftu$ has finite excess equal to $n$. By
changing the index set we may write $\ftu = \fn \cup \fo $, where $\fo$
is a Riesz basis for $\h$. Let $A$ be a lower frame bound for $\fo$ and
choose $\mu < \sqrt{A} $. By compactness there exists a number $m>n$ such 
that $$ || \smo c_i (f_i - g_i ) || \leq \mu \sqrt{ \smo |c_i |^2} $$
for all sets of sequences $ \{c_i \} \subseteq l^2(N) .$ So by 
the remark after Theorem 9,
$\gno$ is a Riesz basis for $\overline{span} \gno .$ If we define the 
operator $T$ on $\h$ by $$Tf_i= f_i ,  \  n<i \leq m , \ \ Tf_i = g_i , 
\ \ i \geq m+1 ,$$ (extended by linearity) then we have an invertible operator
on $\h$. The argument is that every $f \in \h$ has a representation 
$f =  \sum_{i=n+1}^{\infty} c_i f_i $, leading to $$|| (I-T)f ||= 
|| \sum_{i=m+1}^{\infty} c_i (f_i - g_i )|| $$ \ $$ \leq 
\mu \sqrt{ \sum_{i=m+1}^{\infty} |c_i |^2 } \leq \frac{\mu}{\sqrt{A}} 
|| \sum_{i=m+1}^{\infty} c_i f_i ||= \frac{\mu}{\sqrt{A}} \cdot ||f|| .$$
As a consequence, 
$$ codim( \overline{span} \gno) = codim( \overline{span}
 \fno ) = m-n .$$ Take $m-n$ independent elements $\gik$ 
outside $ \overline{span} \gno .$ Then $ \gik \cup \gno$ is a 
frame for $\overline{span} \{ \gik \cup \gno \}
= \h $, since only finitely many elements have been taken away from the
frame $\gtu$. If now $ \sum_{k=1}^{m-n} c_k g_{i_k} + \smo c_i g_i =0 $,
then all coefficients are zero; first,
$$ \sum_{k=1}^{m-n} c_k g_{i_k} = - \smo c_i g_i =0 $$
(if the sums were not equal to zero we could delete an element $g_{i_k}$ and 
still have a frame for
$\h$ contradicting the fact that $ codim( \overline{span} \gno ) =m-n $) and 
since $\gik$ is an independent set and $\gno$ a Riesz basis, all coefficients 
must be zero. So $\gik \cup \gno$ is a Riesz basis, i.e., $\gtu$ also has 
excess $n$.   \\ \\
Now suppose that $\ftu$ has infinite excess. Let $\fti$ be a subset which
is a Riesz basis. Then the corresponding set $ \gti$ spans a space of
finite codimension, i.e., $codim (\overline{span} \gti) < \infty .$
This follows by the same compactness argument as we used in the finite excess
case, which shows that there exist finitely many $f_i , i \in I$  with
the property that if we take them away then we
obtain a family which spans a space with the same codimension 
as the corresponding space of $g_i$'s . Now take a finite family 
$\{g_i \}_{i \in J}$ such that $\{ g_i \}_{i \in I \cup J}$ is total.
Since $\{f_i \}_{i \in I \cup J}$ is a frame with finite excess, the finite
excess result gives that $\{g_i \}_{i \in I \cup J}$ is a frame containing
a Riesz basis, implying that $\gtu$ has infinite excess. \\ \\
We can express the result in the following way: define an equivalence 
relation $\sim $ on the set of frames for $\h$ by
$$ \ftu \sim \gtu \Leftrightarrow K \ \mbox{is compact as an operator 
from} \ l^2(N) \ \mbox{into} \ \h .$$
The equivalence relation partitions the set of frames into equivalence classes.
If a frame contains a Riesz basis, then every frame in its equivalence class 
contains a Riesz basis, and the frames have the same excess. \\ \\
Let us go back to Theorem 10. If $\gtu$ is not total, the rest of the 
assumptions implies that $\gtu$ is a frame for its closed span.Now the
proof of Theorem 10 shows that $\gtu$ contains a Riesz basis for 
$\overline{span} \gtu$, and that the excess refering to this space is
equal to the excess of $\ftu$ as a frame for $\h$ {\it plus} the dimension
of the orthogonal complement of $\overline{span} \gtu$ in $\h$.
\\ \\
Now we want to study the excess property of perturbations in the sense of
Theorem 9 b). We need a result, which might be interesting in itself. 
To motivate it, consider a near-Riesz basis $\ftu$ containing a Riesz basis
$\fti$. Unfortunately, the lower bound for $\fti$ can be arbitrarily small
compared to the lower bound $A$ of $\ftu$. Our result states, that if we are
willing to delete sufficiently (still finitely) many elements, then we
can obtain a family which is a Riesz basis for its closed span, and which
has a lower bound so close to $A$ as we want: \\ \\ 
{\bf Proposition 11:} {\it Let $\ftu$ be a near-Riesz basis with lower bound 
$A$.
Given $\epsilon > 0$ , there exists a finite set $ J \subseteq N $ such that 
$\{f_i \}_{i \in N-J}$ is a Riesz basis for its closed span, with lower
bound $A- \epsilon $.} \\ \\
{\bf Proof:} As in the proof of Theorem 10, write $\ftu = \fn \cup \fo $,
where $\fo$ is a Riesz basis for $\h$. Let $d(\cdot , \cdot )$ denote the
distance inside $\h$ (i.e., $ d(f,E)= inf_{g \in E} ||f-g|| \ \mbox{for} 
\ f \in \h , E \subseteq \h$) and 
choose a number $m> n$ such
that $$ d(f_j , span \{f_i \}_{i=n+1}^m ) < \sqrt{\frac{\epsilon}{n}},
\ \ j=1,...n. $$
We want to show that $\fno$ is a Riesz basis for its closed span, with lower
bound $A- \epsilon .$ Let $P$ dennote the orthogonal projection onto 
$\overline{span} \{ f_i \}_{i=n+1}^m $. Since 
$ || \sum c_i f_i || \geq || \sum c_i (I-P)f_i || $ \ \
for all sequences, it suffices to show that $\{(I-P)f_i \}_{i=m+1}^{\infty}$
satisfies the lower Riesz basis condition with bound $A- \epsilon .$ Let 
$f \in (I-P)\h $. Then 
$$ \smo |< f, (I-P)f_i >|^2 = \siu |<f, (I-P)f_i >|^2 - \sum_{i=1}^n 
|< f, (I-P)f_i >|^2 $$ \ $$ \geq A \nft - 
\sum_{i=1}^n \nft \cdot ||(I-P)f_i ||^2 \geq (A- \epsilon) \nft .$$
Now we only have to show that $ \{ (I-P)f_i \}_{i=m+1}^{\infty}$ is
$\omega$-independent. But if $ \smo c_i (I-P)f_i =0 $, then
$\smo c_i f_i = P \smo c_i f_i $ , implying that both sides are equal to
zero, since $P \smo c_i f_i \in span \{f_i \}_{i=n+1}^m$ and $\fo$ is
independent. Therefore $c_i =0 $ for all $i$. \\ \\    
{\bf Theorem 12:} {\it Let $\ftu$ be a frame for $\h$ with bounds $A,B$. 
Let $\gtu \subseteq \h$ and assume that there exist $\lambda , \mu \geq 0 $
such that $\lambda + \frac{\mu}{\sqrt{A}} <1 $ and   
$$ || \sum c_i (f_i -g_i ) || \leq \lambda \cdot ||  \sum c_i f_i ||
+ \mu \cdot \sqrt{\sum |c_i |^2 } $$ for all finite sequences 
$\{c _i \} $. Then 
$$\ftu \ \mbox{is a near-Riesz basis} \Leftrightarrow 
\gtu \ \mbox{ is a near-Riesz basis}, $$ in which case  
$\ftu$ and $\gtu$ have the same excess.} \\ \\
{\bf Proof:} First assume that $\ftu$ is a near-Riesz basis with
excess $n$. Let $m$ be choosen as in the proof of Proposition 11,
corresponding to an $\epsilon$ satisfying the condition 
$\lambda + \frac{\mu}{\sqrt{A- \epsilon}} < 1$. Let
$Q$ denote the orthogonal projection onto $\overline{span} \fno $. Then 
every element $f \in \h$ can be written $f= (I-Q)f + Qf = (I-Q)f +
\smo c_i f_i $, for some coefficients $c_i$. Now define an operator
$T : \h \rightarrow \h $ by $$ Tf=f , \ f\in \overline{span} {\fno}^{\bot} , 
\ \ Tf_i =g_i , \ i \geq m+1.$$ 
\ $T$ is bounded. Given $f \in \h$ we choose a representation as above. Then
$$ || (I-T)f|| = || \smo c_i (f_i -g_i )|| \leq 
\lambda \cdot || \smo c_i f_i || + \mu \cdot \sqrt{ \smo |c_i |^2} $$ \
$$ \leq  ( \lambda + \frac{\mu}{\sqrt{A- \epsilon}}) || \smo c_i f_i || =
(\lambda + \frac{\mu}{\sqrt{A- \epsilon}}) || Qf|| \leq
(\lambda + \frac{\mu}{\sqrt{A- \epsilon}}) \nf .$$ 
It follows, that $T$ is an isomorphism of $\h$ onto $\h$. So $\gno$ is a
Riesz basis for its closed span, and
$$ dim( \overline{span}{\gno}^{\bot})= dim (\overline{span} {\fno}^{\bot}) .$$
As a consequence, $\ftu$ and $\gtu$ have the same excess. \\ \\
Now assume that $\gtu$ is a near-Riesz basis. By reindexing we may again
assume that $\{g_i \}_{i=n+1}^{\infty}$ is a Riesz basis for $\h$. Define a
bounded operator $U: \h \rightarrow \h $ by $Uf:= \siu <f, \si f_i >g_i .$
Then as in the original proof from \lbrack C3\rbrack \ one proves that
$U$ is an isomorphism of $\h$ onto $\h$ . If we define $U_n : \h \rightarrow
\h $ by $U_n f = \sne <f, \si f_i >g_i $, then this operator has a range
with finite codimension in $\h$, which we will write as
$$ codim_{\h}( R_{U_n}) < \infty .$$
Now let $\etu$ be the natural basis for $l^2 (N)$, i.e., $e_i$ is the 
sequence with $1$ in the i'th entry, otherwise $0$. There exists a bounded
invertible operator $V : \h \rightarrow 
\overline{span} \{e_i \}_{i=n+1}^{\infty} $ such that $Vg_i =e_i $ for 
$i \geq n+1 $, and clearly
$$ codim_{\overline{span} \{e_i \}_{i=n+1}^{\infty}} (R_{VU_n}) < \infty .$$
Observe that $VU_n f = \sne <f, \si f_i >e_i = 
\{ <f, \si f_i > \}_{i=n+1}^{\infty} . $ So
$$ (VU_n)^{*} \{c_i \} = \sne c_i \si f_i = \si \sne c_i f_i .$$
Since $ R_{VU_n}^{\bot} = N_{(VU_n)^{*}}$ has finite dimension, also
$ \{c_i \}_{i=n+1}^{\infty} \longmapsto \sne c_i f_i $ has a finite
dimensional kernel. Therefore 
$$T : l^2(N) \rightarrow \h , \ \ T \ctu = \siu c_i f_i $$
has a finite dimensional kernel, and now the theorem of Holub implies
that $\ftu$ is a near-Riesz basis. By the first part of the Theorem the two 
frames $\ftu$
and $\gtu$ now have the same excess, and the proof is complete. \\ \\ 
Unfortunately, the requirement that $\ftu$ has finite excess is needed in
Theorem 12. In fact we are able to construct examples, where $\{f_i \}$ has
infinite excess and $\{g_i \}$ does not contain a Riesz basis, but where the
perturbation condition is satisfied. Let us shortly
describe how one can do this. Define $\ftni$ as in Lemma 7. 
Given $\epsilon >0$, let
$$g_i^n = e_{\frac{(n-1)n}{2}+i} - \frac{1- \epsilon}{n} 
\sjn e_{\frac{(n-1)n}{2} +j}, \ 1 \leq i\leq n $$ \ 
$$ g_{n+1}^n = \frac{1}{\sqrt{n}} \sjn e_{\frac{(n-1)n}{2}+ j} .$$
Now, given a sequence $ \{\cin \}$ we have 
$$ || \sum \cin (f_i^n - g_i^n ) || =  \epsilon \cdot || \snu [ \sum_{i=1}^n 
\cin] \frac{1}{n} \sju e_{\frac{(n-1)n}{2} +j} || $$ \ $$ \leq    
\epsilon \sqrt{ \snu | \sum_{i=1}^n \cin \frac{1}{\sqrt{n}}|^2 }
\leq \epsilon \sqrt{ \sum | \cin |^2}. \ \ \ \ \ \ \ \ \ \ \ \ \ \ (1) $$
By Lemma 5, $ \{f_i^n \}_{i=1,n=1}^{n, \infty}$ is a frame with bounds $1$.
If we choose $\epsilon <1$, then the perturbation condition is satisfied with 
$\lambda =0, \mu = \epsilon$, implying that 
$ \{g_i^n \}_{i=1, n=1}^{n+1, \infty}$ is a frame with bounds
$(1 - \epsilon)^2 , (1+ \epsilon)^2 .$ \\ \\
{\bf Claim:} $\{g_i^n \}_{i=1, n=1}^{n, \infty}$ is a Riesz basis for $\h$.
\\ \\ We only need to prove that $ \{g_i^n \}_{i=1,n=1}^{n, \infty}$ satisfies
the lower Riesz basis condition.  Given a sequence $\{ \cin \}$ we have
$$ || \snu \sum_{i=1}^n \cin g_i^n || \geq  
|| \snu \sum_{i=1}^n \cin e_{\frac{(n-1)n}{2} +i} || -
(1- \epsilon) || \snu ( \sum_{i=1}^n \cin ) \frac{1}{n} 
\sum_{i=1}^n e_{\frac{(n-1)n}{2}+i} || $$ \ $$ \geq \sqrt{\sum | \cin |^2} -
(1- \epsilon) \sqrt{\snu |\sum_{i=1}^n \cin \frac{1}{\sqrt{n}}|^2} \geq
\epsilon \sqrt{\sum | \cin |^2 } .$$
So actually we have an example where $\ftni$ does not contains a Riesz
basis but the perturbed family does. To obtain the example we were looking
for, we use that $\{g_i^n \}$ has the lower bound $(1- \epsilon)^2 $. 
By (1) above we can consider $\{f_i^n \} $ as a perturbation of $\{g_i^n \}$
if $\frac{\epsilon}{ 1- \epsilon} < 1 $, i.e., if  $\epsilon < \frac{1}{2} $.
So we get our example by choosing $\epsilon < 1 $ and switching the roles of
$\{f_i^n \}$ and $\{g_i^n \}$.  \\ \\
{\it Acknowledgements:} The second author would like to thank Chris Heil for
fruitful discussions on the subject.    
\begin{large}
\begin{center}
References:
\end{center}
\end{large}
\lbrack A\rbrack \ Aldroubi,A.: {\it Portraits of frames.} Proc. Amer. 
Math. Soc., {\bf 123} (1995), p. 1661-1668. \\ \\
\lbrack CC\rbrack \ Casazza, P.G. and Christensen, O.: {\it Hilbert space
frames containing a Riesz basis and Banach spaces which have no subspace
isomorphic to $c_0$.} Submitted, 1995. \\ \\
\lbrack C1\rbrack \ Christensen, O.: {\it Frames and the projection method.}
Appl. Comp. Harm. Anal. {\bf 1} (1993), p.50-53. \\ \\
\lbrack C2\rbrack \ Christensen, O.: {\it Frames containing Riesz bases and
approximation of the frame coefficients using finite dimensional methods.}
Accepted for publication by J. Math. Anal. Appl. \\ \\
\lbrack C3\rbrack \ Christensen, O.: {\it A Paley-Wiener theorem for frames.}
Proc. Amer.Math. Soc. {\bf 123} (1995
), p.2199-2202. \\ \\
\lbrack CH\rbrack \ Christensen, O. and Heil, C.: {\it Perturbation of Banach
frames and atomic decompositions}. Accepted for publication by
Math. Nach.. \\ \\
\lbrack DS\rbrack \ Duffin, R.J. and Schaeffer, A.C.: {\it A class of 
nonharmonic Fourier series.} Trans. Amer. Math. Soc. {\bf 72} (1952)
p. 341-366.  \\ \\
\lbrack HW\rbrack \ Heil, C. and Walnut, D.: {\it Continuous and discrete 
wavelet transforms.} SIAM Review {\bf 31} (1989), p.628-666. \\ \\
\lbrack H\rbrack \ Holub, J.: {\it Pre-frame operators, Besselian frames
and near-Riesz bases.} Proc. Amer. Math. Soc. {\bf 122} (1994), p. 779-
785. \\ \\ 
\lbrack LT\rbrack \ Lindenstrauss, J. and Tzafriri, L.: {\it Classical Banach 
spaces 1.} Springer 1977. \\ \\
\lbrack Se\rbrack \ Seip, K.: {\it On the connection between exponential
bases and certain related sequences in $L^2 (- \pi, \pi )$.} J. Funct. Anal.
{\bf 130} (1995), p.131-160. \\ \\
\lbrack Si\rbrack \ Singer, I.: {\it Bases in Banach spaces 1. } Springer
New York, 1970. \\ \\
\lbrack Y\rbrack \ Young, R.M.: {\it An introduction to nonharmonic Fourier 
series}. Academic Press, New York, 1980. \\ \\
{\bf Peter G. Casazza \\
Department of Mathematics \\
University of Missouri \\
Columbia, Mo 65211 \\
USA. \\
Email: pete@casazza.cs.missouri.edu \\ \\ \\
Ole Christensen \\
Mathematical Institute \\
Building 303 \\
Technical University of Denmark \\
2800 Lyngby \\ 
Denmark. \\
Email: olechr@mat.dtu.dk}

\end{document}